\documentclass{amsart}
\usepackage{amsmath}
\usepackage{amsfonts}
\usepackage{amssymb}
\usepackage{xcolor}
\usepackage{graphicx}
\usepackage[colorlinks=true]{hyperref}
\newtheorem{Theorem}{Theorem}[section]

\newtheorem{Definition}[Theorem]{Definition}
\newtheorem{rem}[Theorem]{Remark}

\newcommand{\R}{\mathbb R}

\newcommand{\N}{\mathbb N}

\newcommand{\F}{\mathcal{F}}

\begin{document}

\title{A gradient method for inconsistency reduction\\
	of pairwise comparisons matrices}

\author{Jean-Pierre Magnot $^{a,b}$}
\author{Ji\v r\'i Mazurek $^{c}$}
\author{Viera \v Cer\v nanov\'a $^{d}$}

\address{$^a$ Univ. Angers, CNRS, \\ LAREMA, SFR MATHSTIC, \\F-49000 Angers, France} 

\address{$^b$ Lyc\'ee Jeanne d'Arc, \\40 Avenue de Grande Bretagne, \\ 63000 Clermont-Ferrand, France}

\address{$^c$ School of Business Administration in Karvina,\\ Department of Informatics and Mathematics,\\ Univerzitni namesti 1934/3, Karvina, Czech Republic}

\address{$^d$ Department of Mathematics and Computer Science,\\ Faculty of Education, Trnava University\\
	Priemyseln\'a 4, 918 43  Trnava, Slovakia}

\email{jean-pierr.magnot@ac-clermont.fr, mazurek@opf.slu.cz, vieracernanova@hotmail.com}

\begin{abstract}
 	We investigate an application of a mathematically robust minimization method - the gradient method - to the consistencization problem of a pairwise comparisons (PC) matrix. Our approach sheds new light on the notion of a priority vector and leads naturally to the definition of instant priority vectors. 
  	
  	We describe a sample family of inconsistency indicators 
  	based on various ways of taking an average value, which extends the inconsistency indicator based on 
  	the ``$\sup$''- norm. We apply this family of inconsistency indicators both for additive and multiplicative PC matrices to show that the choice of various inconsistency indicators lead to non-equivalent consistencization procedures.
\end{abstract}

\maketitle


\section*{Introduction}

	Since human judgments and preferences are often inaccurate and subsequently inconsistent, decision makers have to handle inconsistencies in the evaluation procedure. Many multiple-criteria decision-making (MCDM) problems can be solved efficiently  using pairwise comparisons (PC) via many MCDM methods, such as the analytic hierarchy (network) process AHP(ANP), VIKOR, PROMETHE, ELECTREE or BWM \cite{Brans1985,Figueira2005,Mardani2015,Rezaei2016,Roy1968,Saaty1980}.
	With these methods, one tries to eliminate, or at least, minimize, (cardinal and ordinal) inconsistency in pairwise comparisons. 

	The inconsistency of pairwise comparisons is evaluated by a suitable inconsistency index{, such as } Saaty's consistency index $ CI $ \cite{Saaty1977}, consistency ratio $ CR $ or Koczkodaj's inconsistency index $ KI $ \cite{Koczkodaj1993}. When the inconsistency is ``low'' (Saaty proposed to use the threshold $ CR = 0.10 $), then a priority vector (i.e. vector of weights) is {determined}. However, when the inconsistency of pairwise comparisons is {not reasonably low}, then a decision-maker has two options: to ask an expert to change their judgments, or to find a pairwise comparison matrix that is consistent enough while being {``as close as possible''} to the original matrix. 
	For the latter option, many approaches (algorithms) of inconsistency reduction in pairwise comparison matrices have been proposed in the literature in recent decades.

	{To our best knowledge,} the first algorithm for inconsistency reduction in pairwise comparisons was proposed in \cite{Holsztynski1996} in 1996. The distance-based algorithm searched for the most inconsistent triad (in terms of Koczkodaj's inconsistency index) and replaced it with a consistent one. The same idea was followed in \cite{Koczkodaj2010}. Since then, many other algorithms have been proposed. The objective of each algorithm for inconsistency reduction is to minimize inconsistency under given threshold, or to find a consistent PC matrix ``closest" to a given inconsistent one. Among the proposed algorithms, Abel et al. \cite{Abel2018} gave a multiple-objective optimization algorithm INSITE, 
	Pankratova and Nedashkovskaya \cite{Pankratova2015} proposed a power transformation of an inconsistent pairwise comparison matrix, 
	Boz\'{o}ki et al. \cite{Bozoki2011} proposed an algorithm that finds the minimum number of pairwise comparisons necessary to be altered to achieve desired (lower) inconsistency via a non-linear mixed-integer optimization problem. {Other algorithms can be found in  Li and Ma \cite{LiandMa2007}, 
	Pereira and Costa \cite{Pereira2015,Pereira2018},  Cao et al. \cite{Cao2008}, Zhang et al. \cite{Zhang2014}, Ben\'itez et al. \cite{Benitez2011,Benitez2012,Benitez2014}, Gao et al. \cite{Gao2011} (a similar ant-like algorithm is proposed in  \cite{Girsang2015}), Negahban et al. \cite{Negahban2018}, Boz\'{o}ki et al. \cite{Bozoki2015}. Kulakowski \cite{Kulakowski} introduced a concurrent inconsistency reduction algorithm based on ideas from \cite{Holsztynski1996} and \cite{Xu1999}, which corrects the most inconsistent triad. Similarly, Xu and Xu \cite{Xu2020} adopted a triad approach in inconsistency reduction. 
	Other approaches and algorithms can be found e.g. in  \cite{Ergu2011,MPSK2020,Szybowski2018,Vasconcelos2019}. 


	The aim of the paper is to propose a new inconsistency reduction algorithm, which is based on the selected inconsistency indicator for the evaluation. Indeed, in many approaches, importance has been given to triads with ``big'' inconsistency, but the way of evaluation of the global inconsistency of pairwise comparisons (by the inconsistency indicator) is left as a secondary problem.

	 In a recent work \cite{KMS2019}, the geometric mean method (GMM) and the eigenvalue method (EVM), which provide different priority vectors for inconsistent PC matrices in general, are questioned, but the nature of the inconsistency indicator is ignored. Our approach is reverse in the following sense: whereas inconsistency needs to be minimized and as evaluations are performed by inconsistency indicators, we intend to find  a way (with  changes as small as possible in a PC matrix) to reduce the ``score of inconsistency'' calculated by the inconsistency indicator. If we get a consistent PC matrix, a~family of weights can be generated via the holonomy theorem \cite{Ma2018-1,Ma2018-3} or by any other method. In other words, determining the weights is a secondary problem in our approach, and the inconsistency indicator is understood here as a potential function (depending on the entries of the pairwise comparisons matrix), which has to be minimized (or, ideally, cancelled) first. According to this perspective, gradient  methods are well-behaved for minimizing a potential. Our work improves a basic gradient method for a sample family of inconsistency indicators. Our choice for this sample family is based on the following constraints:
\begin{enumerate}
	\item the inconsistency indicators of this family must be ``similar enough'', in order to ensure that  differences that may occur reveal real sensibility on the choice of the indicator;  
	\item this family must contain at least one well-established inconsistency indicator. 
\end{enumerate}
   	Under these constraints, we choose variations of the well-established Koczkodaj's inconsistency index, which measures the maximal inconsistency of the triads. Our testing family is completed by average inconsistency of triads, average taken in a family that includes arithmetic, quadratic, cubic and harmonic averages. This leads us to make a remark and engage a brief critique on focusing on the only triad with maximal inconsistency, which is given in the appendix. 
   
   	More pragmatically, the direction of the gradient of the inconsistency indicator shows the most efficient way to move the elements of a~PC matrix, in order to achieve a~PC matrix with a~lower value of an inconsistency indicator. This leads us to propose new concepts: an instant priority vector and its discrete analog that we call a difference priority vector. The priority vector shows which state must be attained to reach consistency, while the instant priority vector shows how the elements of a PC matrix must move to get globally a more consistent evaluation.

	Let us finish this description 
 	by a short explanation of the terminology. In our approach, the instant and difference priority vectors provide the desired evolution of pairwise comparisons bringing the PC matrix closer to the set of consistent PC matrices, according to the filtration generated by the inconsistency indicator. Therefore, priority vectors are not mathematically the same objects as instant (or difference) priority vectors, but in applications, our terminology fits with the spirit of (classical) priority vectors: the consistencization of pairwise comparisons.}   

\vskip 12pt
The paper is organized as follows: 
\begin{itemize} 
	\item In section \ref{s:prel}, we describe basics on pairwise comparison matrices, priority vectors and Koczjkodaj's inconsistency index. 
	\item In section \ref{s:lp}, we describe $l^p$ norms, based on the $p-$average of $n$ values, and explain a way of generalizing Koczkodaj's inconsistency index. This generates the announced family of {{}$p-$inconsistency indicators}. 
	\item In section \ref{s:grad}, we sketch a gradient method for a generic piecewise $C^1$ inconsistency indicator and define instant priority vectors, both for additive and multiplicative pairwise comparisons. We focus on our family of {{}$p-$inconsistency indicators}, which helps us to provide calculations in explicit examples of instant priority vectors.
	\item Finally, since it is difficult to determine the integral curves of the gradient flow in general, it is possible to ask for a numerical scheme to obtain an approximate solution: this is the announced gradient method. In addition, the gradient of the inconsistency indicator may be difficult to describe, or in the case of a not-so-regular inconsistency indicator, the gradient may be undefined in classical function theory. We therefore propose to replace the gradient, which is the first derivative, with its discrete analog defined by the difference operator. This leads to the notion of a~difference priority vector, which in most cases is an approximate instant priority vector (section \ref{s:num}). With this last concept, we perform numerical simulations of our method with our sample family  of {{}$p-$inconsistency indicators}, which shows  different results analyzed in the conclusion section. 
\end{itemize}

\section{Preliminaries on pairwise comparisons}\label{s:prel}

Pairwise comparison {{} on  a set of $n$ objects, say indices $i\in \N_n = \{1,\dots,n\},$ is a map from $\N_n^2$ to $\R_+^*=(0,+\infty)$} such that every pair $(i,j)\in \{1,\dots,n\}^2$ is assigned a value $a_{i,j} \in \R_+^* = (0,+\infty), $ {{}which} satisfies the reciprocal condition: {{} $$ \quad a_{i,j} = a_{j,i}^{-1}.$$ These maps are usually represented in a matrix form. Indeed,} 
	the collection of the values $a_{i,j}$ forms an $n \times n$ matrix $A=(a_{i,j})$ called a~(multiplicative) pairwise comparisons matrix, or PC matrix for short.  It is easy to explain the inconsistency in pairwise comparisons when we consider cycles of three comparisons called a ``triad". They are represented here as $(x,y,z)$ with $ x,y,z \in \mathbb{R}_+^* $,
which do not have the ``morphism of groupoid'' property such as $x.z \neq y.$
Therefore, a PC matrix is called consistent if $$\forall (i,j,k)\in \{1,\dots,n\}^3, \quad a_{i,j}a_{j,k} = a_{i,k}, $$ and inconsistent if not. 
	
	\begin{rem}
		Mostly for simplification of calculations to our opinion, some authors also considered additive PC matrices $(b_{i,j})$ obtained from multiplicative PC matrices $(a_{i,j})$ via the group isomorphism $\log : (\R_+^*,\cdot)\rightarrow (\R,+)$ by
		$$\forall (i,j)\in \{1,\dots,n\}^2, \quad b_{i,j} = \log (a_{i,j})\; \Longleftrightarrow \; a_{i,j} = \exp(b_{i,j}).$$ For our preliminary section, we focus on multiplicative PC matrices.
	\end{rem}
A consistent $n\times n$ PC matrix can be derived from a family of weights $(w_i)_{i \in \{1,\dots,n\}}$ called a~priority vector, see e.g. \cite{KMRSTSWY2016}, such that $$\forall (i,j)\in \{1,\dots,n\}^2, \quad a_{i,j} = \frac{w_i}{w_j}.$$ Of course, such a family is not unique (see the obvious correspondence with the mathematical notion of projective coordinates).  For a more theoretical analysis from a different viewpoint, see e.g. \cite{Wajch2019}.

We can create a priority vector even if the PC matrix $(a_{i,j})$ is not consistent, such that   $$\forall (i,j)\in \{1,\dots,n\}^2, \quad a_{i,j} \sim \frac{w_i}{w_j},$$
which is equivalent to finding a consistent PC matrix ``not so far away'' from the original. There are various methods {{} to derive}  a priority vector {{} from} an inconsistent PC matrix. The eigenvalue method and the geometric mean method are the most common ones, their comparison is given in \cite{KMS2019}. 

The apparent fuzziness of our last descriptions shows the need for explicit measurement of the inconsistency of a PC matrix. This is done by selecting the  inconsistency indicator $ii,$ which is a map from the set $PC$ of {{} $\R_+^*-$valued} PC matrices, that  satisfies $$ ii(A) = 0 \quad \Longleftrightarrow \quad A \hbox{ is consistent}.$$  If $ii(A)>ii(B),$ then the PC matrix $A$ is said more inconsistent than the PC matrix $B.$ The family $$\{ii^{-1}([0,\epsilon)) \, | \, \epsilon > 0\}$$ forms a {{} basis of} filter in the set of PC matrices \cite{Ma2018-3}, which hence defines a basis of neighborhoods of the set of consistent PC matrices. The inconsistency indicator can be interpreted as a distance of a given PC matrix to the space $CPC$ of consistent PC matrices.

Koczkodaj's inconsistency index is easily expressed on triads (the upper triangle of a PC submatrix $3 \times 3$):
$$Kii_3(x,y,z) = 1-\min\left \{\frac{y}{xz},\frac{xz}{y}\right \} = \max \left\{ 1- \frac{y}{xz}, 1 - \frac{xz}{y}\right\}.$$

\noindent  According to \cite{KS2014}, it is equivalent to:
$$Kii_3(x,y,z)=1- e^{-\left|\ln\left (\frac{y}{xz}\right )\right |}$$
extending to $n \times n$ PC matrices by the formula:
$${{} Kii_n}(A)=1-\min_{1\le i<j<k\le n}  \min\left(\frac{a_{ik}}{a_{i,j}a_{j,k}}\, ,\,
\frac{a_{i,j}a_{j,k}}{a_{i,k}} \right ) .$$
 
\section{Overview of $l^p-$norms in $\R^n$ and a family of inconsistency indicators}\label{s:lp}
	We consider {{}a family, indexed by $p,$ of} $p-$averages of {{} a fixed vector} $(x_1,\dots ,x_n)\in (\R^*_+)^n, $ where $\R^*_+ = (0,+\infty),$ defined by
$$ || (x_1,\dots,x_n) ||_p = \left( \frac{1}{n} \sum_{i = 1}^n x_i^{p}\right)^{1/p}$$
	for $p \in \R^*.$ 
	Notation $||.||_p$ is justified by the fact that for $p \in [1,+\infty), $ $ ||.||_p$ derives from the well-known $l^p$ norm in $\R^n,$ which is a norm in a standard sense, i.e. it satisfies the triangle inequality. For $p \in (1,+\infty), $
 	the $l^p$-norm is strictly convex. This will be developed thereafter, as this aspect is important for the properties of priority vectors.
  	When $p \in (0;1),$  $||.||_p$ is derived from a~non-convex ``norm", which is a~less developed framework in the mathematical literature. We recall that
\begin{enumerate}
	\item if $p=-1,$ $||(x_1,\dots,x_n)||_p$ is the harmonic mean of $x_1,\dots,x_n,$
	\item if $p=1,$ this is the arithmetic mean,
	\item if $p=2,$ we get the quadratic mean,
	\item  and $||(x_1,\dots,x_n)||_\infty:=\lim_{p \rightarrow +\infty } || (x_1,\dots,x_n) ||_p = \max \{x_1,\dots,x_n\}.$
\end{enumerate} 

	Let us illustrate the geometric differences of the means $|| (|x|,|y|)||_p$ for $p \in \{\frac{1}{2} ; 1; 2; \infty\}$ by level lines. 
	We can see that for $p=1$ or $p=\infty$ the level lines are not strictly convex, and for $p=1/2$ the ball is not convex, {{} see figure 1.} In terms of identities, $||.||_{1/2}$ does not fulfill the triangle inequality, as well as neither of $||.||_p$ for $0<p<1.$ Geometrically, there is no uniqueness of distance projection on convex subsets, while this is true and well-known for $p>1.$ These facts can explain why these ``non convex norms'' (quasi-norms) are less studied.  
\begin{figure}[h!]
	\centering
	\includegraphics[height=2.7cm]{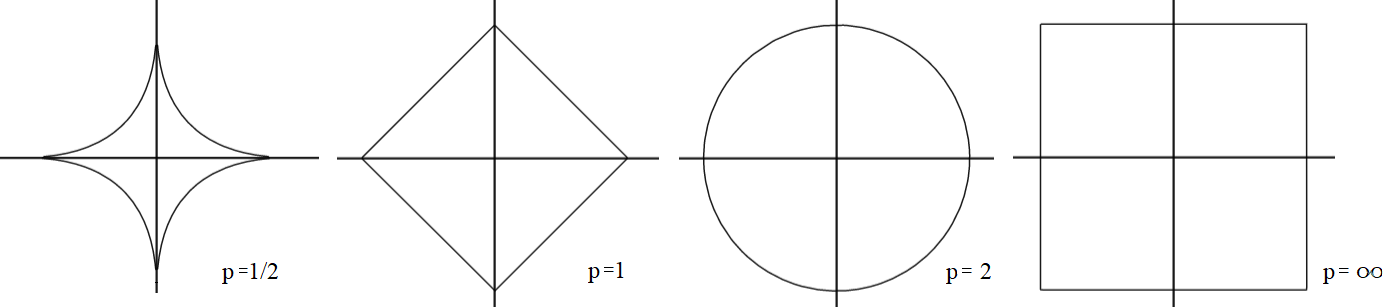}
	\caption{\label{etiquette} Sample level lines for $p=1/2,p=1,p=2$ and $p=\infty.$}
\end{figure}

We generalize Koczkodaj's inconsistency index to the following family of indicators:
\begin{Definition} \label{def:Kii-n-p}
	Let $p \in \R^*.$ We define on $PC_n,$ {{} which is} the set of $n\times n$ pairwise comparisons matrices, the {{}$p-$inconsistency indicator}	
	{{}
	$$ Kii_{n,p}(A) =  1-\exp\left( -\left(\left(\begin{array}{c} n \\ 3 \end{array}\right)^{-1}\sum_{1\leq i < j < k \leq n} \big| \ln(a_{i,j}) + \ln(a_{j,k}) - \ln(a_{i,k})\big|^p\right)^{1/p}\right)  .$$}
	Therewith, we set {$Kii_{n,\infty}= Kii_n$ for the sake of consistent notations.}
\end{Definition}
{{} For this definition, we have to explain that $\left(\begin{array}{c} n \\ 3 \end{array}\right)$ represents the number of possible choices of three indices $i<j<k$ among the $n$ possible indices in $\N_n.$ For example, for $n=4,$ 
	$$\left(\begin{array}{c} n \\ 3 \end{array}\right) = \frac{ 4 \times 3 \times 2}{1 \times 2 \times 3} = 4.$$}
We notice that this inconsistency operator is smooth on an open dense subset of $PC_n,$ and continuous everywhere on $PC_n.$ 

For $p=\frac{1}{2},$ the gradient vectors of $Kii_{n,p}$ will minimize inconsistencies by concentrating them around one axis, i.e. on one triad, which is the reverse effect of the one observed for $p = \infty. $ 

\section{Gradient method for inconsistency indicators and instant priority vector}\label{s:grad}
\subsection{Continuous  and discrete consistencization}
	 Since an inconsistency indicator $ii$ is a real-valued function with a~minimum, a~gradient method (discrete or continuous) can be performed in order to minimize the inconsistency indicator if it is regular enough. For this, the inconsistency indicator needs to be differentiable at least on the set of inconsistent matrices. This method, as a numerical method, can be performed by two ways:
\begin{enumerate}
	\item First, what we call a {{}\textbf{continuous}} gradient method in the multiplicative version of a PC matrix. Here we describe integral curves of the gradient flow in order to get, when intersecting level lines $\{ii(A) = \epsilon\}$ for $\epsilon \geq 0$ and for $A \in PC,$ ``better" PC matrices with ``acceptable" inconsistency ($\epsilon>0$) or consistent PC matrices ($\epsilon = 0$). In case of inability to find an integral curve, we construct an approximate trajectory by generating a~sequence $(A_n)$ of PC matrices by induction. 
	Let us fix $A_0 =(a_{i,j}^{(0)})= A$ and define, for $n \in \N$ and for $i>j,$
	$$ w_{i,j}^{(n)} =  - \frac{\partial}{\partial a_{i,j}} ii(A_{n}),$$
	$$a_{i,j}^{(n+1)} = a_{i,j}^{(n)} + h w_{i,j}^{(n)},$$
	with a constant step $h>0.$ 
	\item Second, a \textbf{discrete} 
	gradient method for the corresponding additive PC matrix $B = (\log(a_{i,j})),$ 
	where we define a sequence $(B_n)$ with a constant step $h>0$ as follows:  
	We start at an additive PC matrix $B_0=B.$ Setting $-w_{i,j}^{(n)}$ as the gradient vector at $B_n,$ we set $$B_{n+1} = \left(b_{i,j}^{(n+1)}\right) \quad \text{and} \quad a_{i,j}^{(n)} = \exp(b_{i,j}^{(n)})$$ with $$ b_{i,j}^{(n+1)}=\log(a_{i,j}^{(n+1)}) = b_{i,j}^{(n)} + h w_{i,j}^{(n)}  ,$$ 
	which is equivalent to 
	{{}$$  a_{i,j}^{(n+1)} = {a_{i,j}^{(n)}}\cdot{\left(\exp(w_{i,j}^{(n)})\right)^h}.  $$}
\end{enumerate}
The stop condition will be of type $ii(A_n) < \epsilon$ with $\epsilon >0,$ because this method only exceptionally provides a consistent PC matrix after a finite number of iterations. 

\begin{rem}
	We have chosen a simple approach with a constant step $ h, $ but for better convergence there are  more sophisticated ways, such as Newton's method.
\end{rem}

\subsection{Gradient and instant priority vector}
 	In decision theory, a priority vector represents the ``best" way to change the coefficients of the PC matrix under consideration in order to make it ``more" consistent, that is, in order to minimize inconsistency. 
 	In classical theory, the priority vector is obtained by a multiplicative mean. 
 	This means that the priority vector does not depend on the selected inconsistency indicator, which encodes whether the inconsistency is acceptable or not.
 The goal of the consistencization procedure is not to get a consistent PC matrix, but to get a PC matrix with lower inconsistency. The minimization of inconsistency should  a~priori depend on the inconsistency indicator under consideration. An effective way to minimize the  inconsistency is to choose the ``best" direction in order to change the elements of the PC matrix. One may ask which ``best" change can be made to achieve maximum change of inconsistency with minimum change in the PC matrix elements. This requirement has a mathematical answer given by the \textbf{gradient}. The gradient of a $C^1-$function $f : \R^n \rightarrow \R$ is defined as 
 $$ \nabla f (x_1,\dots,x_n)= \left(\frac{\partial f}{\partial x_1}(x_1,\dots,x_n),\dots, \frac{\partial f}{\partial x_n}(x_1,\dots,x_n)\right).$$
 \begin{Definition}
 	The \textbf{instant priority vector} with respect to a $C^1-$inconsistency indicator $ii: PC_n \rightarrow \R_+$ is defined by:{{}
 	$$ iW((a_{i,j})_{1\leq i < j \leq n}) = -\nabla ii((a_{i,j})_{1\leq i < j \leq n}).$$}
 \end{Definition}
We then get a $C^0-$vector field over 
$PC_n$ identified with the open domain  $$(\R_+^*)^{\frac{n(n-1)}{2}} \subset  \R^{\frac{n(n-1)}{2}}$$ (we recall that a PC matrix is completely determined by its strict upper triangular part, which represents $1+2+ \cdots + (n-1) = \frac{n(n-1)}{2}$ elements).
\subsection{Instant priority vector on $PC_3:$ case study}
	We give explicit formulas for the announced examples. 
\begin{enumerate}
	\item Let us consider a~generic $3 \times 3$ multiplicative PC matrix
$$A = \left(\begin{array}{ccc} 1 & x & y \\ x^{-1} & 1 & z \\ y^{-1} & z^{-1} & 1 \end{array} \right).$$
	All inconsistency indicators $Kii_{3,p}$ for $p \in \R^*_+$ reduce to $Kii_{3,1}.$
	When $xz/y \neq 1, $ this identification also holds for $p<0.$ 
	The corresponding instant priority vector is, for $xz/y \neq 1, $ 
\begin{equation} 
	iW_{3,1}^m = \text{sign}(y/xz - xz/y) \exp(-|\ln(x)+\ln(z)-\ln(y)|)\cdot \left(-\frac{1}{x} , \frac{1}{y}, -\frac{1}{z}\right),
\end{equation}
	where $\text{sign}(y/xz - xz/y)$ is the sign of $y/xz - xz/y.$ {{}The letter $m$ means that we have a~multiplicative instant priority vector.} Simplifying the notations, 
\begin{equation}
	\left\{	\begin{aligned}
			\hbox{if }\; y/xz > xz/y, \quad &\text{then} \quad iW_{3,1}^m =  \left(-\frac{z}{y} , \frac{xz}{y^2}, -\frac{x}{y}\right),\\
			\hbox{if }\; y/xz < xz/y, \quad &\text{then} \quad iW_{3,1}^m =  \left(\frac{y}{x^2z} , -\frac{1}{xz}, \frac{y}{xz^2}\right).
		\end{aligned}
	\right.	
\end{equation}
	Let us implement the discrete gradient method algorithm with the aim to generate a sequence $(A_n)$ of multiplicative PC matrices
$$A_{n} = \left(\begin{array}{ccc}
	1 & x_n & y_n \\
	x_n^{-1} & 1 & z_n \\
	y_n^{-1} & z_n^{-1} & 1 
\end{array}\right).$$ 
	If $A_0\in PC_3$ is an initial matrix with entries $(x_0,y_0,z_0),$ 
	 we define the three sequences $(x_n),$ $(y_n)$ and $(z_n)$ of positive real numbers by  
$$ \left(
x_{n+1} , y_{n+1} , z_{n+1}\right) =  (x_n,y_n,z_n) + h \cdot iW_{3,1}^m(x_n,y_n,z_n),\quad \forall n \in \N,$$
where $h>0$ is a constant step.


	\item  As a second option, consider an additive PC matrix $$ B = \left(\begin{array}{ccc}
		0 & a & b \\
		-a & 0 & c \\
		-b & -c & 0 \end{array}\right)$$
	with $(a,b,c)=(\log x, \log y , \log z).$ We have 
$$ii(B) = 1 - \exp\left(-|a+c-b|\right)$$ 
and the {\it additive} instant priority vector is
	 $$iW_{3,1}^a = \text{sign}(a+c-b)\exp\left(-|a+c-b|\right)\cdot \left(-1;-1;-1\right).$$
\end{enumerate}
\subsection{From $PC_3$ to $PC_n$}
Let $A=(a_{i,j})$ be a multiplicative $n\times n$ PC matrix. {{}Let $$N =  \left(\begin{array}{c} n \\ 3 \end{array}\right).$$
	If $p \notin \{0;1\},$ the partial derivatives of $||.||_p$ read as:
	$$\frac{\partial}{\partial x_i} \left( \left( \frac{1}{N}
	\left(x_1^p+\dots+x_N^p\right)\right)^{\frac{1}{p}}\right) \; = \;  \frac{1}{N^{1/p}}x_i^{p-1}\cdot\left(\frac{1}{N}\left(x_1^p+\dots+x_N^p\right)\right)^{\frac{1-p}{p}}.$$}
Thus {{}considering $Kii_{n,p}$ as in Definition \ref{def:Kii-n-p},}  
$$\begin{aligned}
	&iW_{n,p}^m(a_{i,j}):=iW_{n,p}^m(A)=\\&= \left(1-Kii_{n,p}(a_{i,j})\right)\cdot\left(-\log(1-Kii_{n,p}(a_{i,j}))\right)^{1-p}\cdot S_{n,p}(A)
\end{aligned}$$ 
	where $S_{n,p}(A)$ stands for  
	$$S_{n,p}(A)={{}\frac{1}{N^\frac{2-p}{p}}}\sum_{i<j<k} \frac{(- \log(1-Kii_{3,1}(a_{i,j},a_{j,k},a_{i,k})))^{p-1}}{1-Kii_{3,1}(a_{i,j},a_{j,k},a_{i,k})}iW^m_{3,1}(a_{i,j},a_{j,k},a_{i,k})$$
	 and $iW^{m}_{3,1}(a_{i,j},a_{j,k},a_{i,k})$ represents  an $\frac{n(n-1)}{2}-$vector with all zeros except in its $(i,j),$ $(j,k)$ and $(i,k)$ components {{}(recall that $p \notin \{0;1\}$)}. 
	 
	 {{}
	 	\begin{rem}
	 		We have also to make a concise precision about the case $p \in \{0;1\}$ and especially $p=1,$ which refers to the use of the arithmetic average. There is a major technical problem: for each $N>2,$ the map $(x_i)_{i \in \N_N}\mapsto ||(x_i)_{i \in \N_N}||_p$ is not differentiable in $\left(\R_+^*\right)^N$ and therefore, one cannot define a (at least) continuous instant priority vector $iW^m_{n,p}$ when $n>3.$ In fact, its definition would be consistent only in a space of generalized functions, noted historically $C^{-1},$ which is part of the space of Schwartz space of distributions. This definition does not seem useful for possible concrete applications of our study.
	 		
	 		This problem will be anyway circumvented in  section 4, since possibly ill-defined derivatives will be replaced by difference operators, which are well-defined for any $p \in \R^*.$ The case $p=0$ seems to us not relevant, because in that case $Kii_{n,0}$ is a constant equal to $1.$ 
	 	\end{rem}}
\subsection{A remark on another possible use of the instant priority vector}
	 Let $CPC_n$ be the set of consistent $n \times n$ PC matrices. The instant priority vectors $iW^a$ and $iW^m$ may be understood as defining a $C^0-$vector field on $(PC_n - CPC_n),$ because $ii$ is a $C^1-$inconsistency indicator.
	 
	 Therefore, one may investigate integral curves for the flow. The flow will define a contraction of $PC_n$ to $CPC_n.$ These integral curves are approximated by the gradient method that we implement. Determining the integral curves is an interesting question. It may enable us to determine a consistencization procedure without the use of algorithms, but that goes far beyond this work. We intent to give new lights on this question in the future.   
\section{Discrete gradient and difference priority vector} \label{s:num}
	A simple way to  numerically implement the gradient method consists in replacing the (continuous) differential by a discrete differentiation, by difference operators. {{} That is, for $n \in \N^*,$  for any differentiable map $f$ from $\R^n$ to $\R,$ and $\forall k  \in \{1,...,n\},$
		\begin{eqnarray*} \frac{\partial f}{\partial x_k}(x_1,...x_n) &=& \lim_{\epsilon \rightarrow 0} \frac{f(x_1,...,x_{k-1},x_k + \epsilon, x_{k+1},...,x_n) - f(x_1,...,x_n)}{\epsilon}\\ &\sim & \frac{f(x_1,...,x_{k-1},x_k + l, x_{k+1},...,x_n) - f(x_1,...,x_n)}{l} \hbox{ for } l \sim 0.\end{eqnarray*}
	The last line defines the so-called difference operator $$ \delta_{x_k}: (f,l) \mapsto \left( (x_1,...x_n) \mapsto \frac{f(x_1,...,x_{k-1},x_k + l, x_{k+1},...,x_n) - f(x_1,...,x_n)}{l}\right)$$  which approximates the classical partial derivative $\frac{\partial}{\partial x_k}$ of each differentiable function $f,$ for each $l$ close enough to $0.$} 
\begin{enumerate}
	\item 
	To approximate {{} the gradient flow, we therefore can replace the classical gradient $\nabla ii$ by a \textbf{discrete} gradient depending on a fixed {{} small number} $l,$ by replacing   $\frac{\partial}{\partial a_{i,j}}ii (A),$ for $A= (a_{i,j})$ by the difference   
\begin{equation} 
	\label{discrete}\delta_{a_{i,j}}(ii)(A,l) = \frac{ii(A') - ii(A)}{l},
\end{equation}
{{}	where} the matrix $A' = (a'_{p,q})$ is defined by :$$a'_{p,q} = a_{p,q}  \quad \text{if} \quad (i,j)\neq (p,q)  \quad \text{and} \quad a'_{i,j} = a_{i,j} + l.$$ }
	Then we can apply the gradient method with respect to this discrete gradient{{}
	$$\delta ii (l) = \left( \delta_{a_{1,1}} (ii,l), ... , \delta_{a_{n,n}} (ii,l)\right)$$ for fixed $l \sim 0.$}
	\item From another viewpoint, one can construct a \textbf{discrete} gradient method on the corresponding additive PC matrix.
\end{enumerate}
	Such methods can be used for {{}$p-$inconsistency indicators} with less regularity, for example when $p=1$ or  $p = \infty.$ If an inconsistency indicator  does not belong to the $C^1-$class, the only problem is the efficiency of the method, while the method  always gives "good" results for $C^1$ inconsistency indicators and $l,$  $h$ ``small enough''.
\subsection{Additive versus multiplicative gradient method on $PC_3$} 
The recursion relation for the iterative method with $iW_{i,j}$ obtained by $(\ref{discrete})$ reads as 
$$ a_{i,j}^{(n+1)} = a_{i,j}^{(n)} + h (iW_{i,j}).$$
Our sample PC matrices are:
	$$ A =  \left(\begin{array}{ccc}
		1 & e^{-2} & e^3 \\
		e^2 & 1 & e \\
		e^{-3} & e^{-1} & 1 \end{array}\right)
	\quad \text{and} \quad 
	B =  \left(\begin{array}{rrr}
		0 & -2 & 3 \\
		2 & 0 & 1 \\
		-3 & -1 & 0 \end{array}\right),$$
	where $\forall (i,j)\in \{1,2,3\}^2, b_{i,j} = \log a_{i,j}.$ The matrix $A$ will test the gradient method on multiplicative PC matrices, while $B$  on additive.  
	In numerical implementation, the iterations become unusable when $ii \sim h$ for $h>>l,$ as show the following tables.
\begin{enumerate}
		\item Results given by the algorithm for the multiplicative matrix $A:$\\
	
\begin{tabular}[b]{|c|c|c|c|}
	\hline
		$h$ & $l$ & best rank $n$ for iteration & ``almost consistent triad'' \\
	\hline
	\hline
		0.1 & 0.001 & $\sim 230$ & 
	\begin{tabular}[c]{@{}l@{}} 
		$a_{1,2}^{(n)} \sim 4,045$ \\ $a_{1,3}^{(n)} \sim 19,676$ \\ $a_{2,3}^{(n)} \sim 4,867$
	\end{tabular}\\
	\hline
		0.01 & 0.001 & $\sim 2300$ & 
	\begin{tabular}[c]{@{}l@{}} 
		$a_{1,2}^{(n)} \sim 4.041$ \\ $a_{1,3}^{(n)} \sim 19.675$  \\ $a_{2,3}^{(n)} \sim 4,868$ 
	\end{tabular} \\ 
	\hline
		0.001 & 0.0001 & $\sim 23000$ & 
	\begin{tabular}[c]{@{}l@{}} 
		$a_{1,2}^{(n)} \sim 4,041$ \\ $a_{1,3}^{(n)} \sim 19,675$ \\ $a_{2,3}^{(n)} \sim 4,868$ 
	\end{tabular} \\
	\hline
\end{tabular}\\
	\item Results given by the algorithm for the additive matrix $B:$\\
	
	\begin{tabular}{|c|c|c|c|}
	\hline
		$h$ & $l$ & best rank $n$ for iteration & ``almost consistent triad'' \\ 
	\hline
	\hline
		0.1   & 0.001  & $\sim 180$ & \begin{tabular}[c]{@{}l@{}}
			$b_{1,2}^{(n)} \sim -0,714$ \\ $b_{1,3}^{(n)} \sim 1,715$ \\ $b_{2,3}^{(n)} \sim 2,285$ 
		\end{tabular} \\ 
	\hline
		0.01  & 0.001  & $\sim 1800$  & \begin{tabular}[c]{@{}l@{}} 
			$b_{1,2}^{(n)} \sim -0,952$ \\ $b_{1,3}^{(n)} \sim 1,120$ \\ $b_{2,3}^{(n)} \sim 2,047$ \end{tabular} \\ 
	\hline
		0.001 & 0.0001 & $\sim 17800$ & \begin{tabular}[c]{@{}l@{}}
			$b_{1,2}^{(n)} \sim -0,667$ \\ $b_{1,3}^{(n)} \sim 1,667$ \\ $b_{2,3}^{(n)} \sim 2,332$ \end{tabular} \\ 
	\hline
\end{tabular}\\
\end{enumerate}

	A simple calculation shows that for the PC matrix after consistencization we are  far from the expected equality $a_{i,j}^{(n)} = \exp(b_{i,j}^{(n)}).$ In addition, it should be noted that
	$$ a_{1,2}^{(n)} < a_{2,3}^{(n)} < a_{1,3}^{(n)},$$
	while 
	$$  b_{1,2}^{(n)} < b_{1,3}^{(n)} < b_{2,3}^{(n)}.$$
	These facts show that the two methods furnish two non-equivalent ways of PC matrix consistencization. Therefore, we now focus on multiplicative PC matrices, moving from $PC_3$ to $PC_4.$ 
\subsection{$p=\infty:$ Koczkodaj's inconsistency index}
\label{infty}
We investigate here the same algorithm for $p=\infty$ and an initial $4 \times 4$ matrix 
$$A = \left(\begin{array}{cccc}
	1& e^{-2} & e^3 & 1 \\
	e^2 & 1 & e & 1 \\
	e^{-3} & e^{-1} & 1 & 1 \\
	1 & 1 & 1 & 1
	\end{array}
	\right)$$

\begin{tabular}{|c|c|c|c|}
	\hline
		$h$ & $l$ & best rank $n$ for iteration & ``almost consistent'' \\ &&& $4\times4$ PC matrix \\
	\hline
	\hline
		0.1 & 0.001 & $\sim 168$ & 
	\begin{tabular}[c]{@{}l@{}} 
		$a_{1,2}^{(n)} \sim 2,517$ \\ $a_{1,3}^{(n)} \sim 19,904$ \\ $a_{2,3}^{(n)} \sim 3,696$ \\ $a_{1,4}^{(n)} \sim 1,398$ \\ 
		$a_{2,4}^{(n)} =1$ \\ 
		$a_{3,4}^{(n)} \sim 0,150$ \\
	\end{tabular}\\
	\hline
	0.01 & 0.001  & $\sim 3080$ & 
	\begin{tabular}[c]{@{}l@{}}
		$a_{1,2}^{(n)} \sim 3,865$ \\
		 $a_{1,3}^{(n)} \sim 19,666$ \\
		 $a_{2,3}^{(n)} \sim 4,812$ \\
		 $a_{1,4}^{(n)} \sim 1,566$ \\
		 $a_{2,4}^{(n)} \sim 0,415$ \\
		 $a_{3,4}^{(n)} \sim 0,083$ 
	\end{tabular}\\
	\hline
\end{tabular}

\subsection{{}$p=1:$ another case of a non-smooth (non-$C^1$) $p-$inconsistency indicator}
With the initial multiplicative PC matrix $A$ from subsection \ref{infty}, we get {{}
$$Kii_{4,1} = 1 - \exp\left(-\frac{1}{4}\sum_{1 \leq i < j<k \leq 4}|\ln(a_{i,j}) + \ln(a_{j,k}) - \ln(a_{i,k})|\right)$$}
and implement the same adapted algorithm.\\ 

\begin{tabular}{|c|c|c|c|}
	\hline
	$h$ & $l$ & best rank $n$  for iteration & ``almost consistent''\\ 	&&& $4\times 4$ PC matrix \\
	\hline
	\hline
	0.1 & 0.001 & $\sim 280	$ & 
	\begin{tabular}[c]{@{}l@{}} 
	 $a_{1,2}^{(n)} \sim 2,768$ \\
	 $a_{1,3}^{(n)} \sim 19,855$ \\
	 $a_{2,3}^{(n)} \sim 3,952$ \\
	 $a_{1,4}^{(n)} \sim 1,544$ \\
	 $a_{2,4}^{(n)} \sim 0,533$ \\
	 $a_{3,4}^{(n)} \sim 0,138$ 
	\end{tabular}\\
	\hline
	0.01 & 0.001  & $\sim 4700$ & 
	\begin{tabular}[c]{@{}l@{}} 
		$a_{1,2}^{(n)} \sim 3,939$ \\
	 $a_{1,3}^{(n)} \sim 19,669$ \\
	 $a_{2,3}^{(n)} \sim 4,812$ \\
	 $a_{1,4}^{(n)} \sim 1,112$ \\
	 $a_{2,4}^{(n)} \sim 0,281$ \\
	 $a_{3,4}^{(n)} \sim 0,057$ 
	\end{tabular}\\
	\hline
\end{tabular}
\subsection{$p \notin \{0;1;\infty\}:$ a smooth case}
We change two last times the {{}$p-$inconsistency indicator}: 
 
 \begin{itemize}
 	\item 
 \underline{$p=2$} uses the quadratic mean {{}
$$ Kii_{4,2}(A) =  1-\exp\left( - \left(\frac{1}{4}\sum_{1\leq i < j < k \leq 4} | \ln(a_{i,j}) + \ln(a_{j,k}) - \ln(a_{i,k})|^2\right)^{1/2}\right)  ,$$}
\begin{tabular}{|c|c|c|c|}
	\hline
	$h$ & $l$ & best rank $n$ for iteration  & ``almost consistent'' \\ &&& $4\times 4$ PC matrix \\
	\hline
	\hline
	0.1 & 0.001 & $\sim 220	$ & 
	\begin{tabular}[c]{@{}l@{}} 
	$a_{1,2}^{(n)} \sim 2,459$ \\
	$a_{1,3}^{(n)} \sim 19,892$ \\
	$a_{2,3}^{(n)} \sim 3,757$ \\
	$a_{1,4}^{(n)} \sim 1,641$ \\
	$a_{2,4}^{(n)} \sim 0,524$ \\
	$a_{3,4}^{(n)} \sim 0,106$ 
	\end{tabular}\\
	\hline
	0.01 & 0.001  & $\sim 3700$ & 
	\begin{tabular}[c]{@{}l@{}} 
	$a_{1,2}^{(n)} \sim 3,571$ \\
	$a_{1,3}^{(n)} \sim 19,725$ \\
	$a_{2,3}^{(n)} \sim 4,573$ \\
	$a_{1,4}^{(n)} \sim 1,613$ \\
	$a_{2,4}^{(n)} \sim 0,422$ \\
	$a_{3,4}^{(n)} \sim 0,089$ 
	\end{tabular}\\
	\hline
\end{tabular}
	\item 
\underline{$p=1/2:$} {{}
$$ Kii_{4,1/2}(A) =  1-\exp\left( - \left(\frac{1}{4}\sum_{1\leq i < j < k \leq 4} | \ln(a_{i,j}) + \ln(a_{j,k}) - \ln(a_{i,k})|^{1/2}\right)^{2}\right)  ,$$}
\begin{tabular}{|c|c|c|c|}
	\hline
	$h$ & $l$ & best rank $n$ for iteration  & ``almost consistent'' \\ &&& $4\times 4$ PC matrix \\
	\hline
	\hline
	0.01 & 0.001 & $=280	$ & 
	\begin{tabular}[c]{@{}l@{}} 
		$a_{1,2}^{(n)} \sim 0,700$ \\
		$a_{1,3}^{(n)} \sim 20,074$ \\
		$a_{2,3}^{(n)} \sim 2,663$ \\
		$a_{1,4}^{(n)} \sim 0,926$ \\
		$a_{2,4}^{(n)} \sim 1,317$ \\
		$a_{3,4}^{(n)} \sim 0,506$ 
	\end{tabular}\\
	\hline
	0.001 & 0.00001  & $\sim 2700$ & 
	\begin{tabular}[c]{@{}l@{}} 
		$a_{1,2}^{(n)} \sim 0,713$ \\
		$a_{1,3}^{(n)} \sim 20,074$ \\
		$a_{2,3}^{(n)} \sim 2,662$ \\
		$a_{1,4}^{(n)} \sim 0,973$ \\
		$a_{2,4}^{(n)} \sim 1,360$ \\
		$a_{3,4}^{(n)} \sim 0,512
		$ 
	\end{tabular}\\
	\hline
\end{tabular}
	\item \underline{$p=-1$} uses the harmonic mean
{{}
	$$ Kii_{4,-1}(A) =  1-\exp\left( - \left(\frac{1}{4}\sum_{1\leq i < j < k \leq 4} | \ln(a_{i,j}) + \ln(a_{j,k}) - \ln(a_{i,k})|^{-1}\right)^{-1}\right)  ,$$}
{ 
	\begin{rem} {{}We notice that  the $p-$inconsistency indicator $Kii_{n,-1}$ is not defined on $CPC_n,$ nor is any $Kii_{n,p}$ indicator  for $p <0.$ 
	Therefore, in contrast to the latest results, we present those for $ h << l, $  because explicit calculations for $ l << h $ fail. This seems to be due to errors caused by the accuracy of the computers. Here, however, we present the problem of mathematical accuracy of the harmonic mean as an open question.
}
\end{rem}  }
\begin{tabular}{|c|c|c|c|}
	\hline
	$h$ & $l$ & best rank $n$ for iteration & ``almost consistent'' \\ 
	&&& $4\times 4$ PC matrix \\
	\hline
	\hline
	0.002 & 0.1 & $=133$ & 
	\begin{tabular}[c]{@{}l@{}} 
	$a_{1,2}^{(n)} \sim 0,228$ \\
	 $a_{1,3}^{(n)} \sim 21,434$ \\
	 $a_{2,3}^{(n)} \sim 2,678$ \\
	 $a_{1,4}^{(n)} \sim 0,991$ \\
	 $a_{2,4}^{(n)} \sim 2,370$ \\
	 $a_{3,4}^{(n)} \sim 0,895$ 
	\end{tabular}\\
	\hline
	0.002 & 0.01  & $=17$ & 
	\begin{tabular}[c]{@{}l@{}} 
	$a_{1,2}^{(n)} \sim 0,144$ \\
	$a_{1,3}^{(n)} \sim 21,737$ \\
	$a_{2,3}^{(n)} \sim 2,713$ \\
	$a_{1,4}^{(n)} \sim 0,999$ \\
	$a_{2,4}^{(n)} \sim 2,654$ \\
	$a_{3,4}^{(n)} \sim 0,986$ 
	\end{tabular}\\
	\hline
\end{tabular}
\end{itemize}

\section{Conclusions}

	By collecting all numerical experiments, it is possible to see that the classical gradient method  used for a selected family of inconsistency indicators:
\begin{itemize}
	\item gives rise to non equivalent consistencization procedures when changing multiplicative to additive PC matrices;
	\item produces different "almost consistent" PC matrices depending on the selected inconsistency indicator,{ see e.g. the sensitivity of the coefficient $a_{2,4}^{(n)}$ for various values of $p,l$ and $h,$ even if our family of test inconsistency indicators are very similar in their way of construction.} 
\end{itemize}
	Therefore, it is necessary to deeply question the classical minimization techniques that do not explicitly depend on the inconsistency: why do they produce the ``best'' and the ``closest'' (almost-)consistent PC matrix? Are they really related to one or the other inconsistency indicator, and if so, in what way? This question is related to the filters around the set of PC matrices, which were described in \cite{Ma2018-3}, where the level sets $ii^{-1}[0;\epsilon)$ form the basis of the filter for the inconsistency indicator $ii$ and $\epsilon >0.$ 

	With the help of the gradient method, we used our selected inconsistency indicators
	like a potential in the space of PC matrices that has to be minimized. The level lines $ii^{-1}(\epsilon)$ are, with the same labeling as before, transverse to the integral curves of the gradient vectors that describe the directions of the ``deepest descent'' between the level lines. Our instant priority vectors naturally depend on the selected inconsistency indicator. {{} This is a striking difference with the work \cite{KSS2020} where orthogonal projections for additive PC matrices are developed, in order to have a global projection scheme from PC to CPC. Although it is well-known that orthogonal projections can be understood as gradient methods for distance functions based on Euclidian norms, the work \cite{KSS2020} appears as incomplete from the viewpoint of our study, since the choice of the orthogonal projection is nowhere better justified than by heuristical arguments of minimization of length. Indeed, a more complete study should have embraced, at least, all the possible positive definite quadratic forms in the space of additive PC matrices, or (even better) the space of Riemannian metrics on the manifold of PC matrices. Even if the authors of this work would have done so, the link with the inconsistency indicator would be still missing in their approach.}  

	This work reinforces the importance of the inconsistency indicator. We hope that it will pave the way to some clarifications on this exact point, which will be remarkable but pragmatic and rigorously well-balanced.

\section*{Appendix: Philosophical aspects of this work: why should we change the paradigm of minimization of maximal inconsistency?} \label{2}
We note that most considered inconsistency indicators express the degree of maximum inconsistency in the pairwise comparison matrix. For example the
Koczkodaj’s inconsistency index can be expressed by the following formula:: 
$$ Kii_3(A) =  1-\exp\left( -\max_{1\leq i < j < k \leq n} | \ln(a_{i,j}) + \ln(a_{j,k}) - \ln(a_{i,k})|\right)  .$$
By the way, in all these choices of inconsistency indicators, minimizing the maximal inconsistency remains the only condition for evaluating the orientations or decisions prescribed, ignoring lower inconsistencies. One can criticize this approach saying that dispatching small inconsistencies in the whole PC matrix is a~better approach than concentrating inconsistency in one triad. In other words, a succession of small liars is better than only one big liar in a speech. Even if psychological issues can justify this choice, 
such an immediate effect cannot be the only justification for not considering all the inconsistencies in the PC matrix. 
For example, trying to find consensual strategies can be the first motivation for trying to make lower kind of average of  the inconsistencies of the PC matrix.   

\paragraph{Acknowledgements} J.M. was  supported by Project no. 21-03085S of the Grant Agency of the Czech Republic. {{} The authors would like to thank the two anonymous referees as well as Eliza Wajch, for their accurate comments which pointed out several inaccuracies in the original version of the text, as well as one of the referees for very interesting remarks about possible developments of this work.}

\end{document}